\documentclass{acm_proc_article-sp}

\usepackage{graphicx}
\usepackage{color}
\usepackage{xypic}

\newcommand{\R}{{\mathbb R}}
\newcommand{\Div}{\text{div}}

\newcommand{\p}{\partial}
\newcommand{\andd}{\quad\text{ and }\quad}
\newcommand{\cancel}[1]{}
\newcommand{\raw}{\rightarrow}

\newcommand{\eps}{\epsilon}

\begin{document}

\conferenceinfo{2009 SIAM/ACM Joint Conference on Geometric and Physical Modeling}{San Francisco, CA}
\CopyrightYear{2009} 

\title{Stable Mesh Decimation}

\numberofauthors{3}

\author{
\alignauthor Chandrajit Bajaj \\
       \affaddr{Department of Computer Sciences}\\
       \affaddr{University of Texas at Austin}\\
       \email{bajaj@cs.utexas.edu}
\alignauthor Andrew Gillette \\
       \affaddr{Department of Mathematics}\\
       \affaddr{University of Texas at Austin}\\
       \email{agillette@math.utexas.edu}
\alignauthor Qin Zhang \\
       \affaddr{Department of Computer Sciences}\\
       \affaddr{University of Texas at Austin}\\
       \email{zqyork@ices.utexas.edu}
}
\date{July 1 2009}

\maketitle
\begin{abstract}
Current mesh reduction techniques, while numerous, all primarily reduce mesh size by successive element deletion (e.g. edge collapses) with the goal of geometric and topological feature preservation.  The choice of geometric error used to guide the reduction process is chosen independent of the function the end user aims to calculate, analyze, or adaptively refine.  In this paper, we argue that such a decoupling of structure from function modeling is often unwise as small changes in geometry may cause large changes in the associated function.  A stable approach to mesh decimation, therefore, ought to be guided primarily by an analysis of functional sensitivity, a property dependent on both the particular application and the equations used for computation (e.g. integrals, derivatives, or integral/partial differential equations).  We present a methodology to elucidate the geometric sensitivity of functionals via two major functional discretization techniques: Galerkin finite element and discrete exterior calculus.  A number of examples are given to illustrate the methodology and provide numerical examples to further substantiate our choices.
\end{abstract}

\section{Introduction}

For function computations carried out on large meshes, mesh decimation is an essential first step.  Mesh decimation techniques are distinguished by the \textit{cost function} they attempt to minimize as they collapse edges in the mesh.  In this paper, we show that given a particular partial differential equation (PDE) problem, an analysis of the geometric sensitivity of the functions involved should guide the choice of cost function for pre-computation decimation.

We consider such function-guided decimation in two realms: adaptive finite element methods (AFEM) and discrete exterior calculus (DEC) methods.  There are three main steps in such methods: formulating a weak version of the governing PDEs, discretizing the problem, and solving the resulting linear system.  Each step introduces a different type of error to the process.  Formulating a weak problem may create what is known as \textit{model error}.  Reducing to a linear system causes \textit{discretization error}.  Implementing the numerical method inevitably causes some  \textit{solver error}.

Adaptive finite element methods aim to control solver error by selective local refinement of the input mesh ($h$-adaptivity), the degree of the basis functions ($p$-adaptivity), or both ($hp$-adaptivity).  While each flavor of adaptive method has met success in particular applications, we will show that their applicability does not immediately transfer to problems that require mesh decimation as a pre-processing step.  Mesh decimation causes a certain loss of geometric information while adaptive refinement is an approximation of missing function information.  For this reason, it is important that the loss of \textit{function} information accrued during mesh decimation be bounded a priori so that the adaptive method can have a hope of converging to a meaningful result.

Discrete exterior calculus methods control solver error by discretizing the functions and operators of the PDE with respect to their algebraic relationships. This type of analysis leads to specific conclusions about where values of the load data and solution data should be assigned or computed; in many cases, values belong most naturally somewhere other than mesh vertices, e.g. on mesh edges or at the circumcenters of triangles.  Therefore, an error bound on function information loss for mesh decimation prior to a DEC method must, by necessity, take into account the locations of the samples of the various variables in the problem.

We describe a framework for selecting an appropriate mesh decimation technique given a PDE and an approach to solving it.  The discretization from the AFEM or DEC method yields a linear system of the form $\mathbb A \textbf x = \textbf b$ whose solution requires inverting the matrix $\mathbb A$.  Our contention is that mesh decimation should be guided by an attempt to avoid large entries in the matrix $\mathbb A$ which can make $\mathbb A$ ill-conditioned and hence destabilize the numerical method.

In Section \ref{sec:pw}, we discuss prior work on AFEM, DEC, and mesh decimation.  In Section \ref{sec:app}, we first give a general overview of AFEM and DEC methods and then explain how each can suggest a mesh decimation technique through a variety of examples.  In Section \ref{sec:res}, we describe the cost functions associated to two existing techniques as well two novel cost functions for use in molecular solvation energetics computations.  In Section \ref{sec:conc}, we present initial experimental results comparing our technique to prior ones.

\section{Prior Work}\label{sec:pw}

We discuss the three main topics of prior work related to our approach: adaptive finite element methods, discrete exterior calculus methods, and mesh decimation methods.

Finite element methods (FEM) have witnessed an explosive growth both in the literature and in industrial application in the past few decades.  Adaptive methods  \cite{BCF1983} have gained traction for their ability to increase local accuracy in a solution.  Beginning with a coarse mesh, AFEM refine by subdividing certain elements into smaller pieces ($h$-adaptivity) \cite{FSPV1989}, increasing the degree of polynomial approximation on some elements ($p$-adaptivity) \cite{BZGd1986}, or a combination of the two ($hp$-adaptivity) \cite{D2007}.  Recently AFEM have been applied to computational biology disciplines.  Baker et al. have worked on a parallel implementation of an AFEM to solve the Poisson-Boltzmann equation \cite{BSHM2000} which generates a mesh of the molecular surface via a subdivision scheme.  Recent work by Chen et al. \cite{CHX2007} provides a FEM for the nonlinear Poisson-Boltzmann equation with rigorous convergence estimates.

Discrete Exterior Calculus (DEC) is an attempt to create from scratch a discrete theory of differential geometry and topology whose definitions and theorems mimic their smooth counterparts.  This theoretical foundation allows the canonical prescription of a discretization scheme for a given PDE problem that enforces topological constraints combinatorially instead of numerically, thereby providing for increased robustness in implementation.    This approach has been employed by an increasing number of authors in recent years to develop multigrid solvers \cite{wB2008}, solve Darcy flow problems \cite{HNC2008}, and geometrize elasticity \cite{Y2008}.  For a complete introduction to DEC theory, see Hirani \cite{H2003} and Desbrun et al. \cite{DHLM2005}.  In this paper, we give a brief introduction to the theory in Section \ref{sec:app} and an example in Section \ref{subsec:darcy}.

Mesh decimation techniques are abundant in geometry processing literature.  A useful survey of many methods was given by Heckbert and Garland \cite{HG1995} and a more recent book by Luebke et al. \cite{LWCRV2002} provides a thorough overview of mesh simplification techniques.   Bajaj and Schikore \cite{BS1998} have given an error bounded mesh decimation technique for 2D scalar field data.  We focus on mesh decimation of unstructured surface meshes via edge contraction including the approaches of Garland and Heckbert \cite{GH1997} and Lindstrom and Turk \cite{LT1998,LT1999}.  These are explained in Section \ref{sec:res}.

\section{Methodology for Elucidating Geometry Sensitive Functionals}
\label{sec:app}
We begin with a general problem: find $u\in V$ such that
\begin{equation}\label{eqn:strong}
Lu = f\quad\text{ on $\Omega$},
\end{equation}
where $L$ is a linear operator, $V$ is the appropriate solution space for the problem, and $\Omega$ is a simplicial complex embedded in $\R^3$.  There are two main techniques used to discretize this into a linear system: Galerkin methods and Discrete Exterior Calculus (DEC) methods.  We describe how each could be used and how a mesh decimation technique should be chosen accordingly.

The Galerkin finite element method begins by putting (\ref{eqn:strong}) into the weak form: find $u\in V$ such that
\begin{equation}\label{eqn:weak}
a(u,v)=(f,v)_{L^2} \quad\forall v\in V,
\end{equation}
where $a$ is the operator $L$ phrased as a bilinear form (usually symmetric) and $f$ is treated as a functional $(f,\cdot)_{L^2}$ on $V$.  An appropriate finite dimensional subspace $V_h\subset V$ is chosen and an answer to the following discretized problem is sought: find $u\in V_h$ such that
\[ a(u,v)=(f,v)_{L^2} \quad\forall v\in V_h. \]
Since $V_h$ is finite-dimensional, we can fix a basis $\{\phi_i:1\leq i\leq n\}$ of $V_h$.  The size of the basis is proportional to the number of elements in the mesh $\Omega$.  Write $u=\sum_{j=1}^nU_j\phi_j$, $K_{ij}=a(\phi_j,\phi_i)$ and $F_i=(f,\phi_i)$.  Set $\textbf{U}=(U_j)$, $\mathbb K=(K_{ij})$ $\textbf{F}=(F_i)$.  Then solving (\ref{eqn:weak}) over $V_h$ is the same as solving the matrix equation
\begin{equation}\label{eqn:femlinsys}
\mathbb K\textbf U = \textbf F.
\end{equation}
For a proof and detailed discussion, see \cite{BS2002}.

A significant amount of care goes into the selection of $V_h$ to ensure that the method is both well-posed and stable.  ``Well-posed'' means the system has a unique solution and ``stable'' means there exists a constant $C>0$ independent of $h$ such that
\[ ||u-u_h||_V \leq C\inf_{w_h\in V_h}||u-w_h||_V. \]
In other words, a method is stable if the error between the true solution $u$ and approximate solution $u_h$ is bounded above uniformly by a constant multiple of the minimal approximation error for $V_h$.  The famous Babuska inf-sup condition \cite{BA1972} is often used to simultaneously prove both well-posedness and stability of a FEM and hence provide an a priori bound on solver error.  We note, however, that the stability error bound does not account for error due to naive mesh decimation.

If the mesh $\Omega$ has too many elements, (\ref{eqn:femlinsys}) will be too large for the solver, making decimation necessary.  For decimation to be useful, however, it must not create very large or small entries in $\mathbb{K}$ which might make $\mathbb K$ ill-conditioned.  Since the entries of $\mathbb K$ are a functional $a(\cdot,\cdot)$ on the basis functions $\phi_i$, mesh decimation must be guided by the geometry-sensitive components of $a$ as opposed to the geometry of $\Omega$ alone.  Such components are necessarily problem-specific as we elucidate in examples presented in Sections \ref{subsec:pbe} and \ref{subsec:gbe}.

The method of Discrete Exterior Calculus is an alternative approach which focuses on correctly discretizing the operator $L$ instead of the solution space $V$.  The viewpoint provided by differential geometry and topology reveals how this ought to be done.  Common operators such as grad, curl, and div are all manifestations of the exterior derivative operator $d$ in dimensions 1, 2, and 3, respectively.  Equations relating quantities of complementary dimensions, such as the constitutive relations in Maxwell's equations, involve a Hodge star operator $\ast$ which provides the canonical mapping.  The Laplacian operator $\Delta$ can be written as $\delta d+d\delta$ where $\delta$ is the coderivative operator, defined by $\delta:=\ast d\ast$.  Each operator has a discrete version designed to mimic the properties of its smooth counterpart.  The discrete versions of the operators are written as matrices whose entries depend only on the topology and geometry of the mesh $\Omega$.

To solve (\ref{eqn:strong}), an analysis is made as to the dimension of $u$ as a $k$-form based on either the problem context or the type of operator acting on it.  The div operator in 3D, for example, acts on 2-forms while grad acts on 0-forms.  The variable $u$ is replaced by a vector $\vec u$ with one entry for each $k$-simplex in the mesh of $\Omega$ and the operator $L$ is replaced by its discrete counterpart, written as a matrix $\mathbb L$.  The load data $f$ is converted to a vector $\vec f$ accordingly.  This yields the equation
\begin{equation}\label{eqn:declinsys}
\mathbb L\vec u = \vec f,
\end{equation}
which can then be solved by linear methods.

Again, it may be necessary to decimate $\Omega$ so that the linear system (\ref{eqn:declinsys}) is tractable on a computer.  The size of the entries of $\mathbb L$ depend heavily on the geometry-sensitive operators such as $\ast$ and $\delta$ and less on the topology sensitive operators such as $d$.  Therefore, to prevent an ill-conditioned $\mathbb L$, mesh decimation must be guided based on the definition of the discrete operators as opposed to the definition of the solution space.  We discuss an example in Section \ref{subsec:darcy}.

\subsection{Poisson-Boltzmann Electrostatics}
\label{subsec:pbe}
The Poisson-Boltzmann equation (PBE) describes the attraction between solvated molecules.  We describe its linearized version according to the formulation given by Lu, Zhang, and McCammon in \cite{LZM2005}, which is believed to be a sufficient approximation for electrostatics computations.  Let $\Omega\subset\R^3$ be a domain indicating interior molecular regions with point charges $q_1,\ldots, q_N$ located at $r_1,$$\ldots,$$r_N\in\Omega$.  The linear PBE is
\begin{equation}
\label{eqn:phiint}
\nabla^2\phi^{\text{int}}(r_p)=-\frac 1{\eps_{\text{int}}}\sum_{k=1}^N q_k\delta(r_p-r_k),\quad p\in\Omega,
\end{equation}
\begin{equation}
\label{eqn:phiext}
\nabla^2\phi^{\text{ext}}(r_p)=\kappa^2\phi^{\text{ext}}(r_p),\quad p\in{\Omega},
\end{equation}
where $\phi^{\text{int}}$ and $\phi^{\text{ext}}$ are the electrostatic potentials on the interior and exterior of $\Omega$, $\eps_{\text{int}}$ is the interior dielectric constant, $\delta$ is the Dirac distribution, and $\kappa$ is the inverse of the Debye-H\"{u}ckel screening length.  Values for the constants $\eps_{\text{int}}$ and $\kappa$ are determined experimentally.   At the boundary $\partial\Omega$, the surface potential $f$ should satisfy $f=\phi^{\text{int}}=\phi^{\text{ext}}$ with normal derivative $h=\p\phi^{\text{ext}}/\p n$.  The equations are put into integral form and discretized, reducing the problem to a set of linear equations of the form
\begin{equation}
\label{eqn:pbematrix}
\left(\begin{array}{ll}\mathbb{B} & \mathbb{A} \\\mathbb{D} & \mathbb{C}\end{array}\right)\left(\begin{array}{l} f \\ h\end{array}\right) = \left(\begin{array}{l} Q \\ 0\end{array}\right)
\end{equation}
where $Q$ is the initial data of point charges and locations.  The four (sub)matrices $\mathbb{A}$, $\mathbb{B}$, $\mathbb{C}$ and $\mathbb{D}$ have entries
\[\sum_{t}\int_{E_t} \mathcal{I} dA\]
where $E_t$ is a facet of $\Omega$, $t$ indexes over a small neighborhood of mesh elements, and the integrand $\mathcal{I}$ is one of the following:
\[\mathcal{I}\in\left\{G(x_i,x_j), \frac{\p G(x_i,x_j)}{\p n}, u(x_i,x_j), \frac{\p u(x_i,x_j)}{\p n}\right\}.\]
The functions $G$ and $u$ are the Green functions for (\ref{eqn:phiint}) and (\ref{eqn:phiext}), respectively.  They are given by
\[G(x_i,x_j)=\frac 1{4\pi r_{ij}}\andd u(x_i,x_j)=\frac{\exp(-\kappa r_{ij})}{4\pi r_{ij}},\]
where $r_{ij}=|x_i-x_j|$.  Hence, the terms of the submatrices in (\ref{eqn:pbematrix}) decay like $1/r_{ij}$ at worst.  To avoid a blowup of these terms, we use a cost function $f_{pb}$ that attempts to preserve mesh element quality as a way of avoiding small $r_{ij}$ values.  We describe $f_{pb}$ in Section \ref{subsec:fpb}.

\subsection{Generalized Born Electrostatics}
\label{subsec:gbe}
A recent approach by Bajaj and Zhao \cite{BZ2009} computes molecular solvation energetics and forces by using a Generalized Born (GB) model instead of a Poisson Boltzmann model of electrostatic solvation.  While the PB model begins with the PBE, a description of the electrostatic potential over the whole domain, the GB model begins with a model of the solvation energy of a single atom in a given medium.  The electrostatic solvation energy of a molecule is defined in terms of the pairwise interaction between these atomic energies:
\[ G_{\text{pol}} = -\frac\tau 2\sum_{i,j}\frac{q_iq_j}{[r_{ij}^2+R_iR_j\exp(-r_{ij}^2/4R_iR_j)]^{1/2}}.\]
Here, $\tau=\frac 1{\eps_p}-\frac 1{\eps_w}$ where $\eps_p$ and $\eps_w$ are the solute and solvent dielectric constants, $q_i$ and $R_i$ are the charge and effective Born radius of atom $i$, and $r_{ij}$ is the distance between atom $i$ and atom $j$.  The success of a GB method hinges upon efficient and accurate approximation of the effective Born radii $R_i$.  Bajaj and Zhao use the surface integration technique given in \cite{GRF1998}.  This gives the expression
\[ R_i^{-1}=\frac 1{4\pi}\int_{\Gamma}\frac{(\textbf r - \textbf{x}_i)\cdot\textbf n(\textbf r)}{|\textbf r-\textbf{x}_i|^4}dS,\quad i=1,\ldots,M,\]
where $\Gamma$ is the solvent-molecular interface, $\textbf{x}_i$ is the center of atom $i$, and \textbf n(\textbf r) is the unit normal of the surface at \textbf{r}.  The position vector \textbf{r} ranges over $\Gamma$.  The integral is approximated by
\[ R_i^{-1}\approx \frac 1{4\pi}\sum_{k=1}^N w_k\frac{(\textbf{r}_k - \textbf{x}_i)\cdot\textbf n(\textbf{r}_k)}{|\textbf{r}_k-\textbf{x}_i|^4}dS,\quad i=1,\ldots,M, \]
where the $\textbf{r}_k$ are the Gaussian quadrature nodes with weights $w_k$ lying on a triangular mesh approximating the surface $\Gamma$.  Therefore, the computation of $R^{-1}_i$ will be sensitive to changes in the position of Gauss points relative to the nearest atomic centers, i.e. changes in the computed values of $|\textbf{r}-\textbf{x}_i|$.  Accordingly, we design a cost function $f_{gb}$ to minimize the cumulative change in $|\textbf{r}-\textbf{x}_i|$ values.  This function is described in Section \ref{subsec:fgb} and compared experimentally to other cost functions in Section \ref{sec:conc}.

\subsection{Darcy Flow}
\label{subsec:darcy}
Recent work by Hirani et al. \cite{HNC2008} uses a DEC method to model Darcy flow, a description of the flow of a viscous fluid in a permeable medium.  The governing equations under the assumption of no external body force are given by
\[
\begin{array}{rcll}
f + \frac k\mu\nabla p & = & 0 & \text{in $\Omega$,} \\
\Div f & = & \phi & \text{in $\Omega$,} \\
f & = & \psi & \text{on $\p\Omega$,}
\end{array}
\]
where $f$ is the volumetric flux, $k>0$ is the coefficient of permeability, $\mu>0$ is the coefficient of viscosity, $\phi:\Omega\raw\R$ is the prescribed divergence of velocity, and $\psi:\p\Omega\raw\R$ is the prescribed normal component of the velocity across the boundary.  They discretize the operators $\nabla$ and $\Div$ based on DEC theory to arrive at the linear system
\[
\left[\begin{array}{cc}
-(\mu/k)\mathbb{M}_{n-1} & \mathbb{D}^T_{n-1} \\
\mathbb{D}_{n-1} & 0
\end{array}
\right]
\left[\begin{array}{c} f \\ p \end{array}
\right]
=
\left[\begin{array}{c} 0 \\ \phi \end{array}
\right].
\]
Here, $\mathbb D_k$ is the discrete exterior derivative operator that acts on $k$-cochains and $\mathbb{M}_k$ is a diagonal matrix representing the Hodge Star operator on $k$-cochains.

We now consider the effect of a priori mesh decimation for this scheme.  The DEC analysis used to derive this method requires that the solution $[f\enspace p]^T$ provide values of the flux $f$ on $(n-1)$-simplicies (i.e. edges in triangle meshes and triangles in tetrahedral meshes) and values of the pressure $p$ at the circumcenters of $n$-simplicies.

For the pressure values to have any meaning, the mesh must be \textit{well-centered} meaning the circumcenter of each simplex must lie in the interior of the simplex.  Since this criterion is often violated by meshing schemes (e.g. an obtuse triangle is not well-centered), pressure is assigned instead to the barycenters of $n$-simplicies.  The authors point out that the error introduced by this modification prevents the exact representation of linear variation of pressure over the domain.  Therefore, an appropriate cost function for this method should be weighted to favor the creation of simplified meshes with good quality elements (e.g. elements with good aspect ratios).  This would minimize the distance between the barycenter and circumcenter, making the calculations more robust.

\section{Description of Cost Functions}
\label{sec:res}

\begin{figure}
\centering
\includegraphics[width=.45\textwidth]{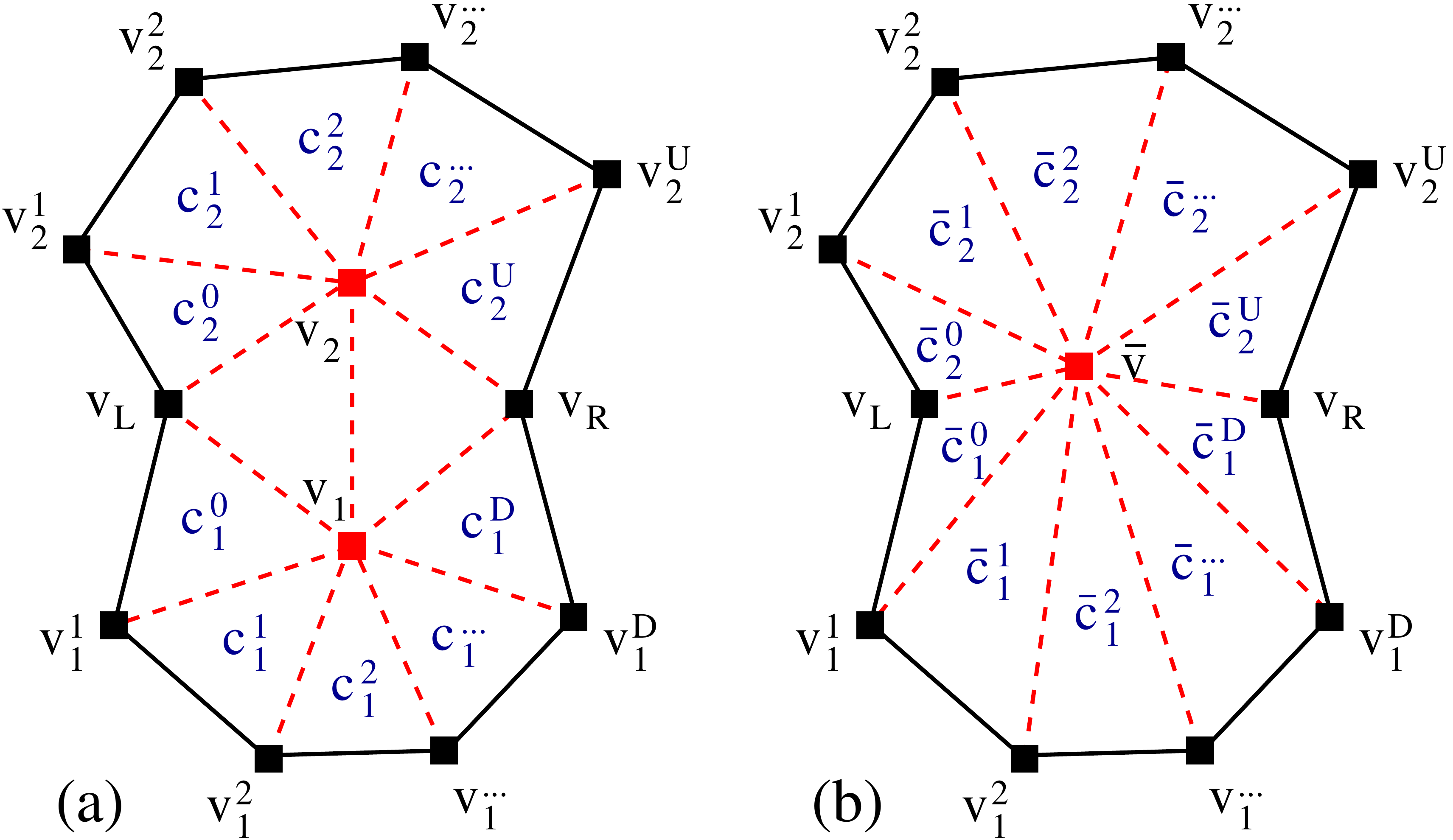}
\caption{Notation for collapse of edge $(v_1,v_2)$ to $\overline{\textbf{v}}$.  Black vertices and edges are unchanged in the collapse while red vertices and dashed edges may change.}
\label{fig:edgeContrac}
\end{figure}
\subsection{Quadratic Error Cost Function}
\label{subsec:fqe}

The quadric error measure we use comes from Garland and Heckbert \cite{GH1997}.  First, an error metric $\Delta(\textbf v)$ is established for each vertex \textbf{v}, based on the planes $P(\textbf v)$ that contain the triangles incident to \textbf{v}.  A plane $\textbf p\in P(\textbf v)$ given by $ax+by+cz+d=0$ is represented as $[a\enspace b\enspace c\enspace d]^T$.  The vertex \textbf{v} is represented as $[v_x\enspace v_y\enspace v_z\enspace 1]^T$.  Then the error metric $\Delta(\textbf v)$ is defined by
\[ \Delta(v) = \textbf{v}^T\left(\sum_{\textbf p\in P(\textbf v)}\textbf p\textbf{p}^T\right)\textbf v. \]
For a vertex $\textbf{v}_i$, let $\mathbb Q_i=\sum_{\textbf p\in P(\textbf{v}_i)}\textbf p\textbf{p}^T$.  Then the cost of collapsing edge $(\textbf{v}_1,\textbf{v}_2)$ to some point $\overline{\textbf{v}}$ is
\[ f_{qe}(\textbf{v}_1,\textbf{v}_2;\overline{\textbf{v}}) := \overline{\textbf{v}}^T(\mathbb{Q}_1+\mathbb{Q}_2)\overline{\textbf{v}}. \]
We use the publicly available software QSlim \cite{qslim} to implement this cost function.  For a given edge, the program attempts to find an optimal placement of $\overline{\textbf{v}}$ by solving a certain linear system derived from the $\mathbb{Q}_i$.  If the matrix associated to this system is not invertible, it tries to place $\overline{\textbf{v}}$ optimally on $(\textbf{v}_1,\textbf{v}_2)$.  If this fails, it sets $\overline{\textbf{v}}$ to be either $\textbf{v}_1$, $\textbf{v}_2$, or the midpoint of the edge, whichever minimizes $f_{qe}$.

\subsection{Volumetric Error Cost Function}
\label{subsec:fvol}

The volumetric error measure we use comes from Lindstrom and Turk \cite{LT1998,LT1999}.  First, we establish the notation for the signed volume $V$ of a tetrahedron bounded by a vertex $\textbf v$ and the vertices $\textbf{v}^{t_i}_0$, $\textbf{v}^{t_i}_1$, $\textbf{v}^{t_i}_2$ of a triangle $t_i$.  As before, vertices are written as four component vectors, e.g. \textbf{v} is represented by $[v_x\enspace v_y\enspace v_z\enspace 1]^T$.  Then $V$ is defined by
\[ V(\textbf{v},\textbf{v}^{t_i}_0, \textbf{v}^{t_i}_1, \textbf{v}^{t_i}_2) = \frac 16\det(\textbf{v}\enspace \textbf{v}^{t_i}_0\enspace \textbf{v}^{t_i}_1\enspace \textbf{v}^{t_i}_2)  =: \frac 16 \textbf{G}_{t_i}\textbf v \]
where $\textbf{G}_{t_i}$ is a 1 by 4 matrix defined by the above equation.  The cost of collapsing edge $(\textbf{v}_1,\textbf{v}_2)$ to some point $\overline{\textbf{v}}$ is
\[ f_{vol}(\textbf{v}_1,\textbf{v}_2;\overline{\textbf{v}}) := \frac 12\overline{\textbf{v}}^T\left(\frac{1}{18}\sum_i \textbf{G}_{t_i}^T\textbf{G}_{t_i}\right)\overline{\textbf{v}},\]
where $i$ indexes over triangles $t_i$ incident to at least one of $\{\textbf{v}_1, \textbf{v}_2\}$.  This cost function is ultimately quite similar to $f_{qe}$ except that $f_{qe}$ weights the distance between $\overline{\textbf{v}}$ and a plane by the area of the triangle defining the plane while $f_{vol}$ weights it by the \textit{square} of the triangle area.  As is explained in \cite{LT1998,LT1999}, this weighting better serves the goal of volume preservation.  We use a package from the publicly available software TeraScale Browser \cite{tsb} to implement this cost function.

\subsection{Poisson Boltzmann Cost Function}
\label{subsec:fpb}
In Section \ref{subsec:pbe} we discuss how Poisson Boltzmann (PB) computations are sensitive to Gaussian quadrature points coming into close proximity.  Hence, we define a cost function $f_{pb}$ which penalizes such occurrences.

We fix the following notation for the collapse of edge $(\textbf{v}_1,\textbf{v}_2)$ to the point $\overline{\textbf{v}}$.  Consider the union of triangles incident to $\textbf{v}_1$ or $\textbf{v}_2$.  These are the only vertices, edges, and triangles whose geometry may be changed by the edge collapse.  Although all these objects lie in $\R^3$, their generic connectivity information is captured in $\R^2$ by Figure \ref{fig:edgeContrac} (a).  Taking $(\textbf{v}_1,\textbf{v}_2)$ to be vertical with $\textbf{v}_1$ at the bottom, the triangle to the left (resp. right) of the edge has $\textbf{v}_L$ ($\textbf{v}_R$) as its third vertex.  We proceed from $\textbf{v}_L$ to $\textbf{v}_R$ along the upper (resp. lower) vertices labeling them $\textbf{v}_2^1$, $\textbf{v}_2^2$, $\ldots$, $\textbf{v}_2^U$ ($\textbf{v}_1^1$, $\textbf{v}_1^2$, $\ldots$, $\textbf{v}_1^D$ (D for down)).  Set $\textbf{v}_2^0:=\textbf{v}_1^0:=\textbf{v}_L$ and $\textbf{v}_2^{U+1}:=\textbf{v}_1^{D+1}:=\textbf{v}_R$.  We denote the centers of the triangles not adjacent to $(\textbf{v}_1,\textbf{v}_2)$ as
\[c_2^u:=\frac 13\left(\textbf{v}_2+\textbf{v}_2^u+\textbf{v}_2^{u+1}\right),\quad u=0,1,\ldots,U,\]
\[c_1^d:=\frac 13\left(\textbf{v}_1+\textbf{v}_1^d+\textbf{v}_1^{d+1}\right),\quad d=0,1,\ldots,D.\]
The collapse operation moves $\textbf{v}_1$ and $\textbf{v}_2$ to $\overline{\textbf{v}}$.  The result is shown in Figure \ref{fig:edgeContrac} (b).  The $\{\textbf{v}_1^d\}$ and $\{\textbf{v}_2^u\}$ are unchanged, but the new centers are $\{\bar c_2^u\}$ and $\{\bar c_1^d\}$ where $\overline{\textbf{v}}$ replaces $\textbf{v}_2$ or $\textbf{v}_1$ in the expression of the center.  Re-index $c_1^d$ and $c_2^u$ as $c_i$.

If the triangles indexed by the $c_i$ are of good quality, meaning nearer-to-equilateral, their Gaussian quadrature points will be better spaced.  We approximate the quality of a triangle $T$ by
\[q(T) := \frac{l(T)}{s(T)}+\frac{maxa(T)}{mina(T)},\]
where $l(T)$ (resp. $s(T)$) denotes the longest (shortest) side and $maxa(T)$ (resp. $mina(T)$) denotes its maximum (minumum) angle.  The best quality triangles have the minimum $q$ value of 2.  We denote triangles in Figure \ref{fig:edgeContrac} by their center $c_i$ or $\bar c_i$.  The cost function is then defined to be
\[f_{pb}(\textbf{v}_1,\textbf{v}_2;\overline{\textbf{v}}) := \sum_iq(T_{\bar c_i})-q(T_{c_i}).\]
To further improve element quality, we decimate in stages and run a quality improvement code based on geometric flow \cite{XPB2006} in between stages.

\subsection{Generalized Born Cost Function}
\label{subsec:fgb}
In Section \ref{subsec:gbe} we discuss how the Generalized Born (GB) computations are sensitive to changes in the location of Gaussian quadrature points $\{c_i\}$ relative to the atomic centers $\{x_j\}$ of the molecule in question.  Hence, we define a cost function $f_{gb}$ which penalizes edge collapses with a larger cumulative change in $|c_i-x_j|$ values.  Since it would be too computationally expensive to compute the complete change for every edge collapse, we use a restricted set of pertinent $\{c_i\}$ and $\{x_j\}$ values described below.

We take the $\{c_i\}$ described in Section \ref{subsec:fpb} as our set of pertinent Gaussian quadrature points.  We set $\{x_j\}$ to be those atomic centers lying within a fixed distance $\rho$ of either $v_1$ or $v_2$.  Since Born radii are on the order of 1-2 \AA, $\rho$ should be set between 2 and 5 to capture a manageable, non-empty set of nearby atoms.  In the future, we will devise a parameter sweep to optimize the value of $\rho$.

We want to minimize the atomic center functional $f_{ac}$
\[f_{ac}:=\left|\sum_{i,j}|c_i-x_j|^2-|\bar c_i-x_j|^2\right|,\]
which can be re-written as
\[f_{ac}=\left|\sum_{i,j}c_i^Tc_i-\bar c_i^T\bar c_i+2(\bar c_i^T-c_i^T)x_j \right|.\]
All the variables in the above expression are known, save for the $\bar c_i$ which are linear functions of $\overline{\textbf{v}}$.  We define the GB-dependent cost function to be
\[f_{gb}(\textbf{v}_1,\textbf{v}_2;\overline{\textbf{v}}) := |\textbf{v}_1-\textbf{v}_2| + \lambda f_{ac}(\textbf{v}_1,\textbf{v}_2;\overline{\textbf{v}}).\]
The first term is used to promote the collapse of shorter edges and thereby improve triangle quality.  The weight factor $\lambda$ is chosen so that the two terms are of the same order of magnitude.

\section{Experimental Results and Conclusions}
\label{sec:conc}
\begin{figure}
\centering
\includegraphics[width=.4\textwidth]{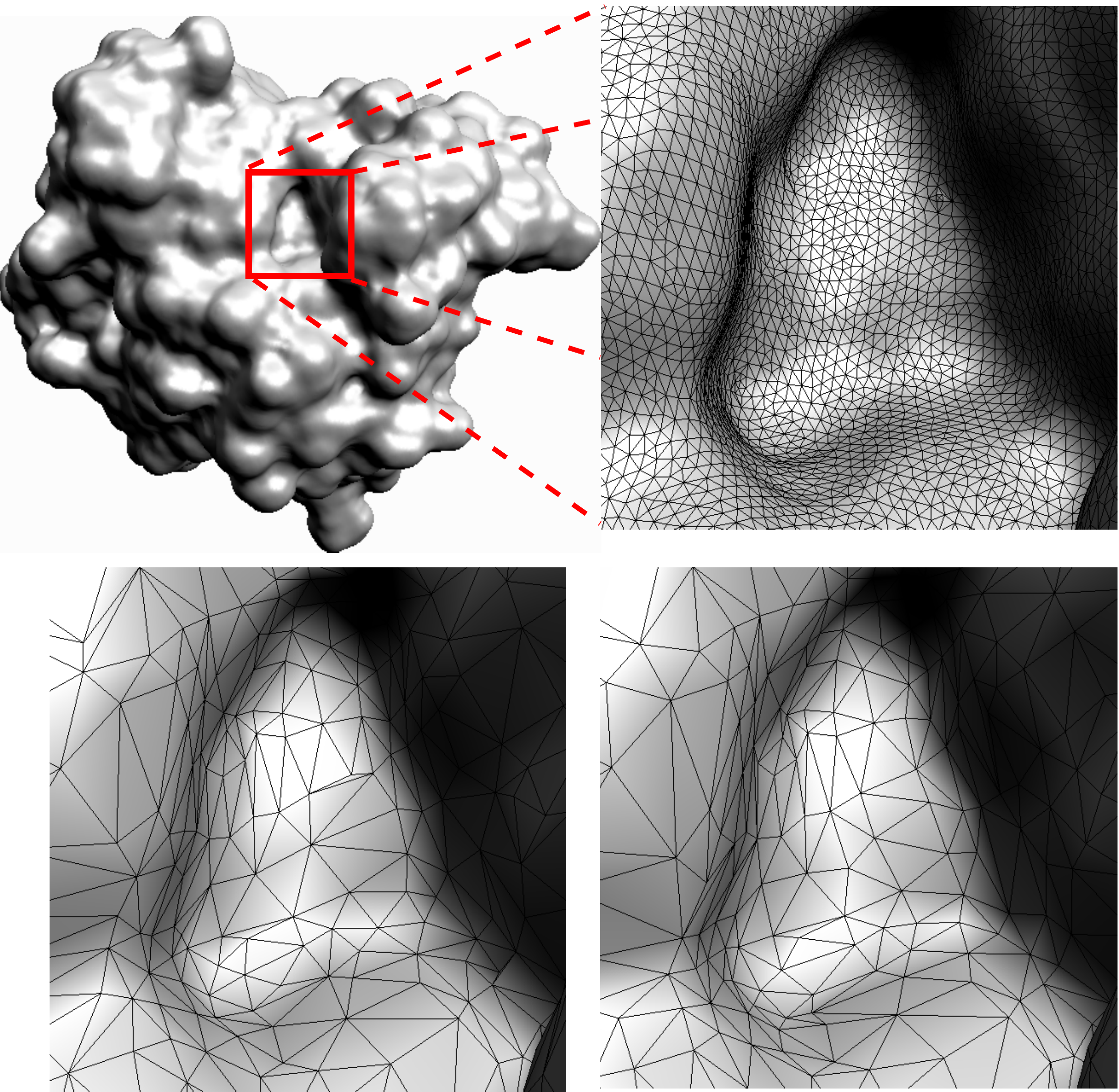}
\caption{\textbf{\textit{Top:}} Surface rendering of mAChE before decimation.  The original mesh of the surface has about 200,000 vertices and 400,000 faces.  In the inset, the fine mesh is visible.  This pocket region of the molecule aids in its biological function. {\textbf{\textit{Bottom:}}} The mesh of the pocket after decimation to 25,500 faces using $f_{qe}$ (left) and a modified version of $f_{gb}$ (right).}
\label{fig:decim_zoom}
\end{figure}
\begin{figure}
\centering
\includegraphics[width=.45\textwidth]{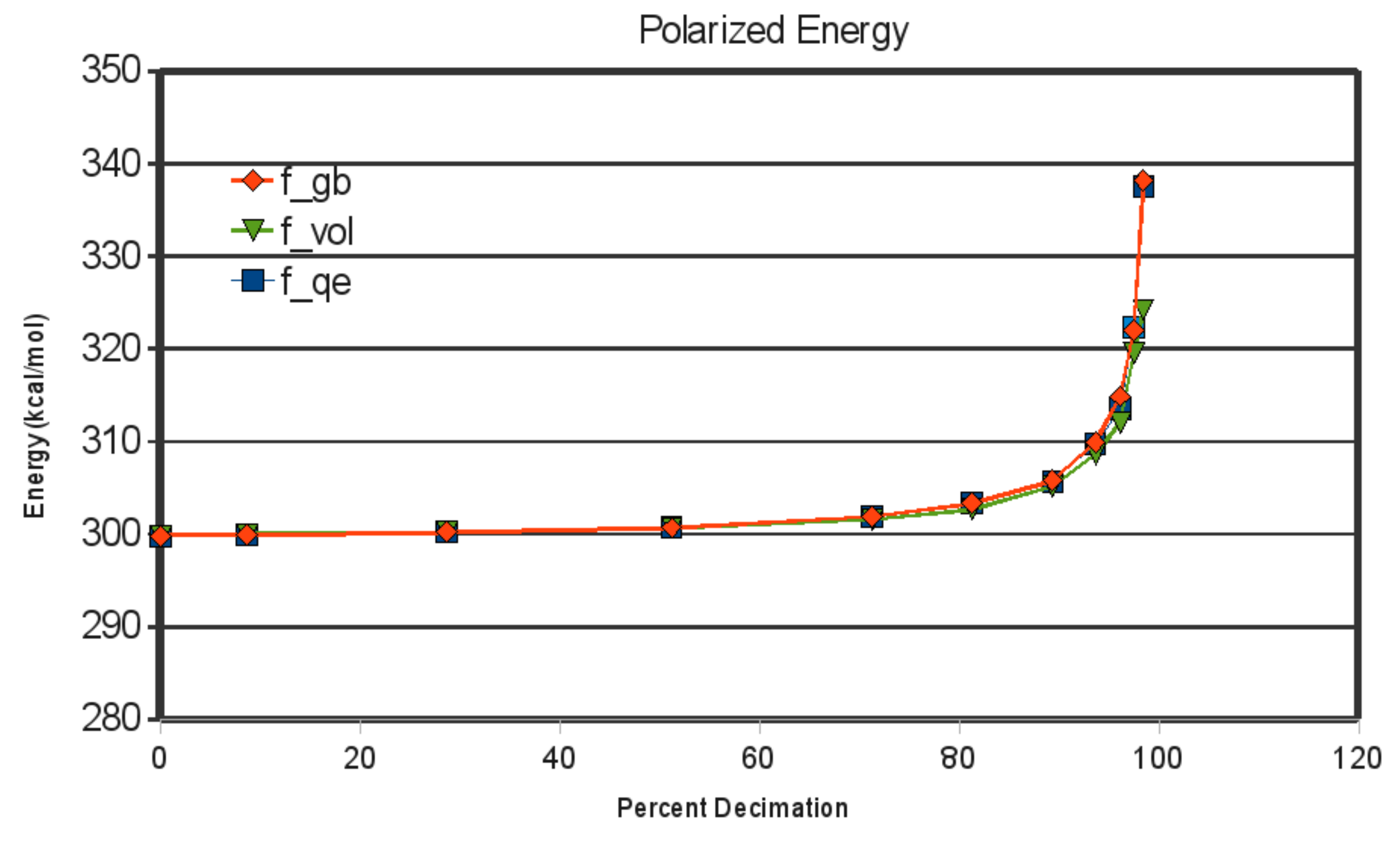}
\includegraphics[width=.45\textwidth]{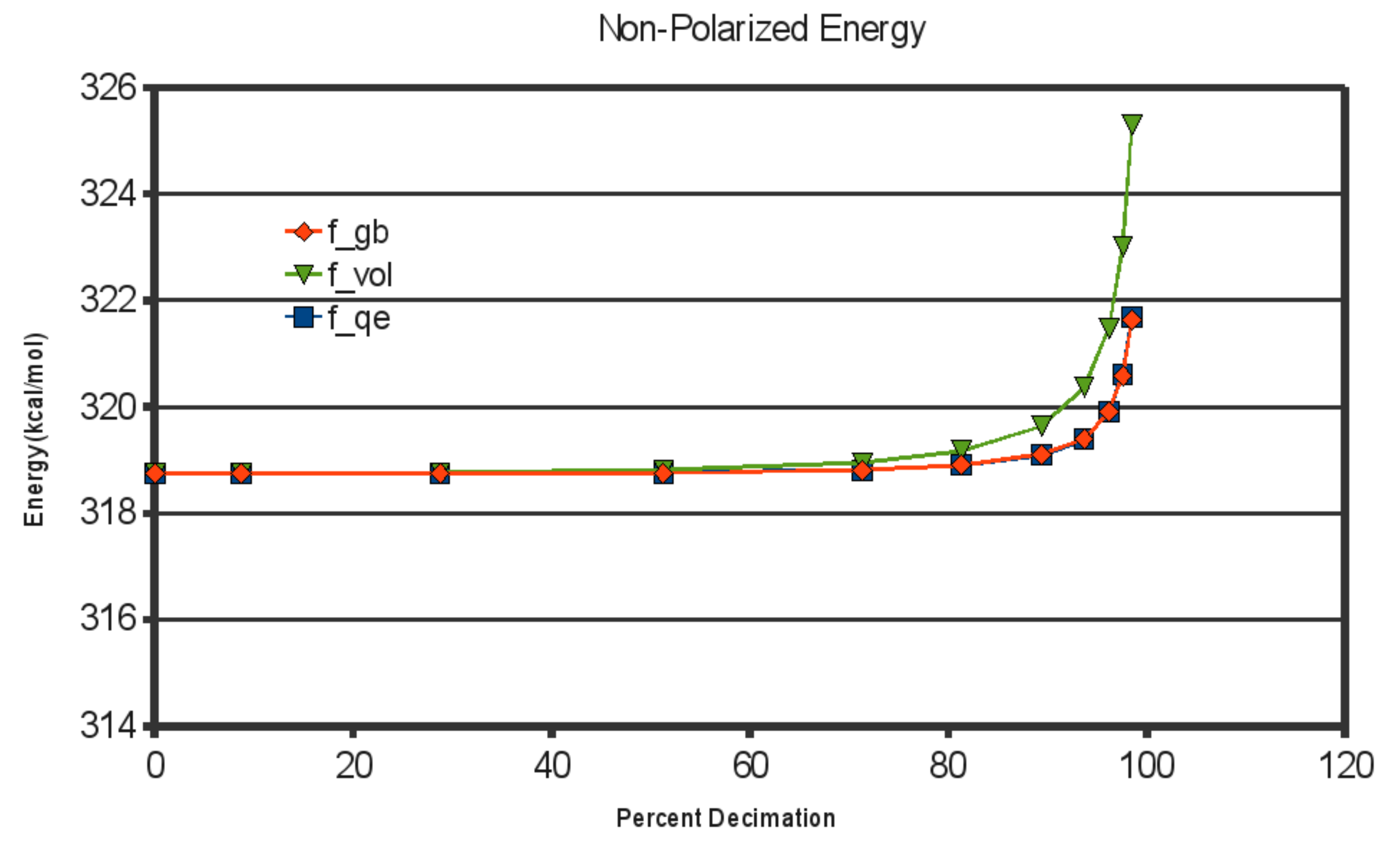}
\caption{Effect of cost function used for decimation on computed GB energy values for mAChE.}
\label{fig:res}
\end{figure}

To test the validity of our claims, we work with Mouse Acetylcholinesterase (mAChE).  This macromolecule serves an important regulatory function as it terminates the action of the neurotransmitter acetylcholine (ACh).  The initial mesh of the molecular surface is generated from a Protein Data Bank (PDB) \cite{pdb} file of the molecule using our in-house software TexMol \cite{texmol}.  We show a picture of the initial mesh in Figure \ref{fig:decim_zoom}.  We decimate this mesh using the cost function $f_{qe}$ and a modified version of $f_{gb}$ which uses $f_{qe}$ instead of $|\textbf{v}_1-\textbf{v}_2|$.  We use a Dynamic Packing Grid data structure \cite{BCR2009} to efficiently compute the set $\{x_j\}$ of nearby centers for each mesh edge.  We set $\rho=5$ and $\lambda=10^{-8}$ and compute the polarized and non-polarized energy for each mesh using the nFFGB code described in \cite{BZ2009}.  The meshes are only marginally different as shown in Figure \ref{fig:decim_zoom} and thus produce similar energy values as shown in the charts in Figure \ref{fig:res}.  With further experimentation and parameter sweeps, we believe $f_{gb}$ will begin to out-perform $f_{qe}$.  Still, Figure \ref{fig:res} shows that $f_{vol}$ is a decidedly worse choice for non-polarized energy computations as it does not respect the functional sensitivity of the problem.  In future work, we will also implement $f_{pb}$ and use PB-CFM code by Bajaj and Chen \cite{BC2009} to compute and compare PB energetics.

\section{Acknowledgments}

We would like to thank Dr. Peter Lindstrom for his help with the TeraScale Browser software and Dr. Wenqi Zhao for her help with the energetics calculations.   This research was supported in part by NSF grants DMS-0636643, CNS-0540033 and NIH contracts R01-EB00487, R01-GM074258, R01-GM07308.

\bibliographystyle{abbrv}
\bibliography{siamacm2009}

\end{document}